\documentclass[12pt]{article}
\usepackage[left=3cm,top=3.0cm,right=3cm,bottom=2.6cm]{geometry}

\usepackage[ansinew]{inputenc}
\usepackage{amscd,amsfonts,amsmath,amssymb,amsthm,enumerate,epsfig,float,graphics,graphicx,pst-all,yhmath}
\newtheorem{lemma}{Lemma}

\newtheorem{theorem}{Theorem}
\newtheorem{corollary}{Corollary}
\newtheorem{example}{Example}
\newtheorem{remark}{Remark}
\newtheorem{definition}{Definition}

\newcommand{\fin}{\hfill $\Box$}

\title {Convex bodies with equipotential circles}
\author{
Iv\'an Gonz\'alez-Garc\'ia$^1$, Jes\'us Jer\'onimo-Castro$^{2}$ \\ Valent\'in Jim\'enez-Desantiago$^{3}$, and Efr\'en Morales-Amaya$^{4}$\\ 
\small{$^{1,2}$Facultad de Ingenier\'ia}\\
\small{Universidad Aut\'onoma de Quer\'etaro, M\'exico}\\
\small{$^{3}$Instituto de Matem\'aticas,}\\
\small{Universidad Nacional Aut\'onoma de M\'exico, M\'exico}\\
\small{$^{4}$Facultad de Matem\'aticas-Acapulco,}\\
\small{Universidad Aut\'onoma de Guerrero, M\'exico}\\
}

\begin{document}

\maketitle
\begin{abstract}  
Given a convex body $K\subset \mathbb R^2$ we say that a circle $\Omega\subset \text{int} K$ is an \emph{equipotential circle} if every tangent line of $\Omega$ cuts a chord $AB$ in $K$ such that for the contact point $P=\Omega\cap AB$ it holds that $|AP|\cdot|PB|=\lambda$, for a suitable  constant number $\lambda$. The main result in this article is the following: \emph{Let $K\subset\mathbb R^2$ be a convex body which has an equipotential circle $\mathcal B$ with centre $O$ in its interior. Then $K$ has centre of symmetry at $O$, moreover, if none chord of $K$ which is tangent to $\mathcal B$ subtends an angle $\pi/2$ from $O$, then $K$ is a disc.} We also derive some results which characterizes the ellipsoid and the sphere in $\mathbb R^3$ and introduce also the concept of equireciprocal disc.
\end{abstract}

\section{Introduction}
\begin{figure}[H]
    \centering
    \includegraphics[width=.42\textwidth]{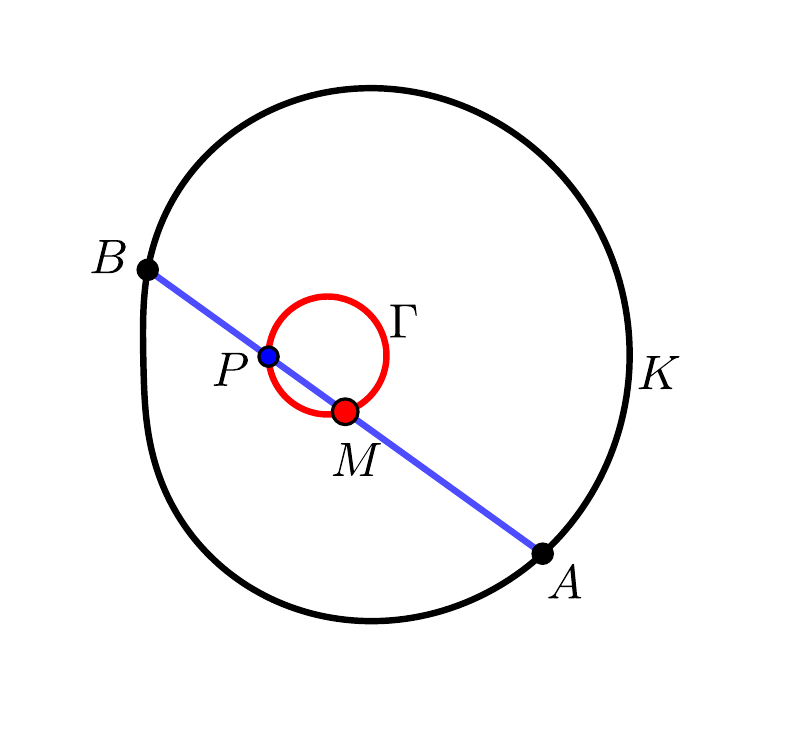}
    \caption{The Lima\c{c}on of Pascal has an equichordal point}
    \label{cardioide}
\end{figure}

Let $\Gamma$ be a closed convex curve in the plane, i.e., the boundary of a compact and convex figure in the plane, and let $P$ be a point in its interior. We say that $P$ is an \emph{equichordal point} if every chord of $K$ through $P$ have the same length. There is a very famous problem due to W. Blaschke, W. Rothe, and R. Weitzenb\"ock \cite{Blaschke}, which asks if there exists a convex body with two equichordal points. There is a long story, full of many false proofs about the no existence of such a body. In \cite{Rychlik} M. Rychlik finally gave a complete proof about the no existence of a body with two equichordal points. However, there are many convex bodies, different from the disc, which have exactly one equichordal point. One example of these bodies is the Lima\c{c}on of Pascal shown in Fig. \ref{cardioide}, moreover, this figure has also the property that the locus of midpoints of the chords through the equichordal point is a circle. 

Another interesting kind of points is the following: we say that $P$ is an \emph{equipotential point} if for every chord $AB$ of $K$, through $P$, the product of the lengths $|AP|\cdot|PB|$ is equal to a given constant value. Clearly, every point in the plane is an equipotential point if $\Gamma$ is a circle. Indeed, if with respect to a given curve $\Gamma$, every point $P$ in the plane is equipotential, then $\Gamma$ is a circle. The proof of this fact is easy and we let it to the interested reader. A more interesting question is: how many equipotential points with respect to a curve $\Gamma$, ensure us that it is  a circle?

The answer to the above question is 2, as was proved by K. Yanagihara in \cite{Ya2}. Some years later, J. Rosenbaum also proved that two equipotential points are enough \cite{Ros}. That only one equipotential point is not enough is shown by the curve already given by Yanagihara \cite{Ya2}. Let $ABCDEF$ be a regular hexagon with centre $O$; consider the circumscribed circles to the triangles $\triangle OBC$, $\triangle ODE$, and $\triangle OFA$. The curve formed by the segments $AB$, $CD$, $EF$, and the arcs $\widehat{BC}$, $\widehat{DE}$, and $\widehat{FA}$, has $O$ as an equipotential point.

\begin{figure}[H]
    \centering
    \includegraphics[width=.42\textwidth]{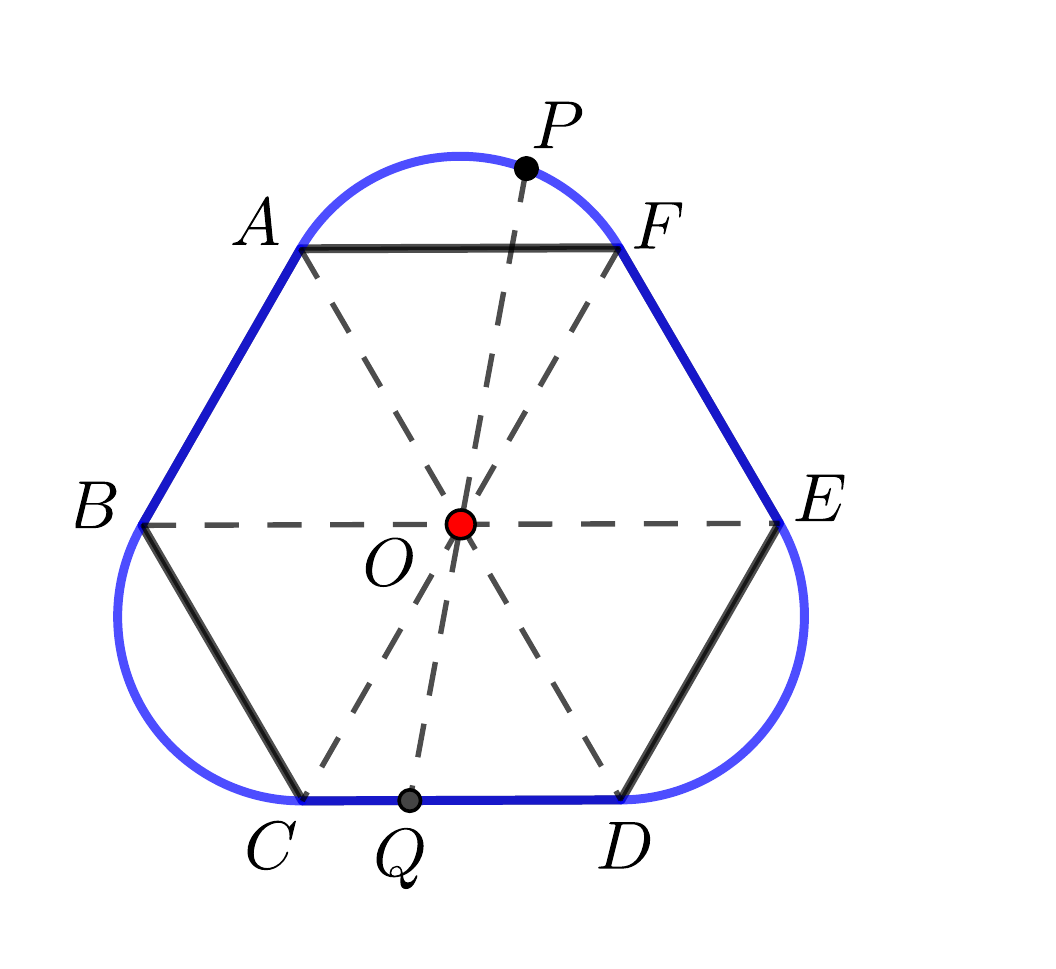}
    \caption{A curve with an equipotential point $O$}
    \label{yanagihara}
\end{figure}

We would like to extend the notion of equipotential and equichordal points to fatter sets of points.

\begin{definition}
Given a convex body $K\subset \mathbb R^2$ we say that a circle $\Omega\subset \emph{int} K$ is an \emph{equipotential circle} (\emph{equichordal circle}) if every tangent line of $\Omega$ cuts a chord $AB$ in $K$ such that for the contact point $P=\Omega\cap AB$ it holds that $|AP|\cdot|PB|=\lambda$ (it holds that $|AB|=\lambda$), for a suitable  constant number $\lambda$.
\end{definition}

In \cite{Barker}, J. A. Barker and D. G. Larman proved that the only convex body in the plane which possesses an equichordal circle is the disc. The question now is, what convex bodies, besides the discs, possesses an equipotential circle?

By the Example 37 of Chapter XII in the book by C. Smith, we know that the envelope of the chords in an ellipse, which are seen under an angle of $\pi/2$ from the centre of the ellipse, is a circle. Indeed, this is the unique equipotential circle for the ellipse. For a given ellipse $K$, this equipotential circle is constructed as follows: suppose the centre of $K$ is $O$, the minor and mayor axes $AC$ and $BD$, respectively. Let $\Omega$ be the circle inscribed in the rombus $ABCD$ and let $P$ be the point of tangency between $AB$ and $\Omega$. Since $\measuredangle AOB=\pi/2$, we have that $|AP|\cdot |PB|=|PO|^2=r^2,$ where $r$ is the radius of $\Omega$. We reformulate our question in a more precise way: what convex bodies, besides discs and ellipses, possesses an equipotential circle?

\begin{figure}[H]
    \centering
    \includegraphics[width=.5\textwidth]{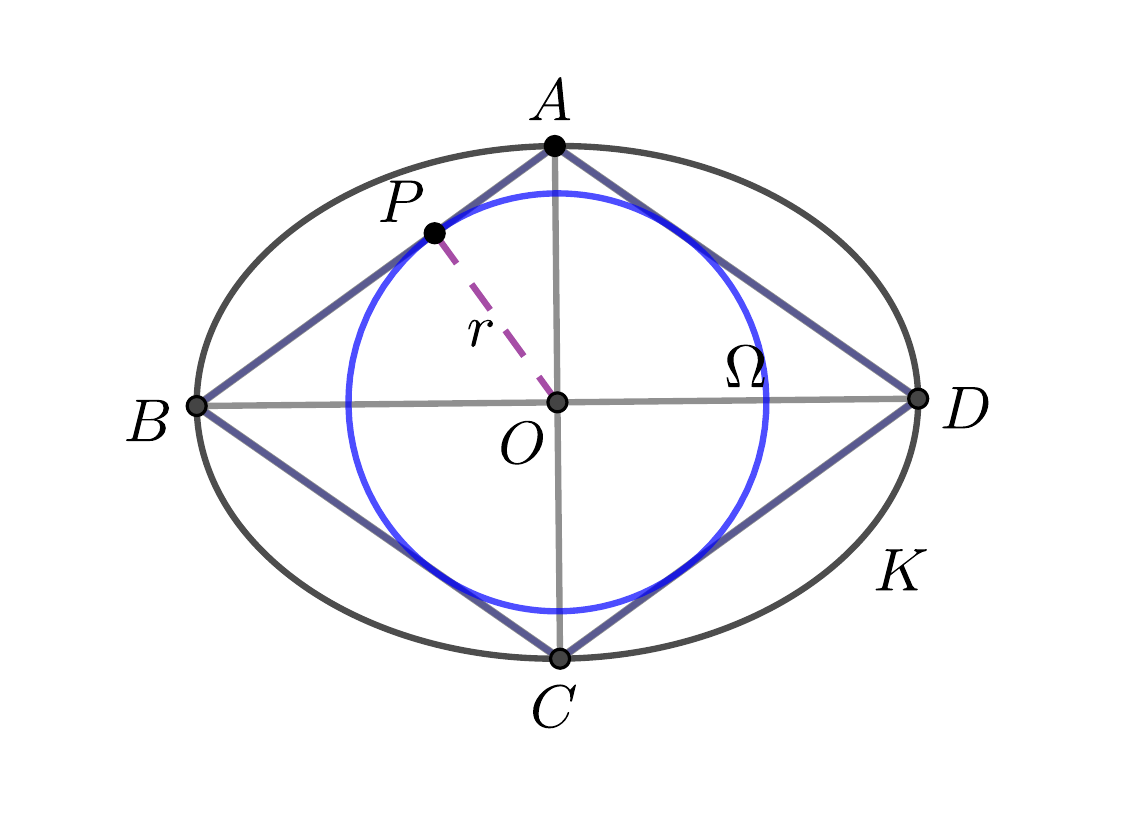}
    \caption{Ellipses have equipotential circles}
    \label{ejemplo_elipse}
\end{figure}

We solved almost completely this question in Theorem \ref{principal} and use the result to give some characterizations of the circle in $\mathbb R^2$, the sphere and the ellipsoid in $\mathbb R^3$. We also deal with one more kind of circles: \emph{equireciprocal circles}.

\section{Convex bodies with equipotential circles}
The main result in this section is the following.

\begin{theorem}\label{principal}
Let $K\subset\mathbb R^2$ be a convex body which has an equipotential circle $\mathcal B$ with centre $O$ in its interior. Then $K$ has centre of symmetry at $O$, moreover, if none chord of $K$ which is tangent to $\mathcal B$ subtends an angle $\pi/2$ from $O$, then $K$ is a disc.
\end{theorem}

Before we give the proof of Theorem \ref{principal}, we will give some useful lemmas. Let $AB$ be any chord of $K$ which is tangent to $\mathcal B$ at a point $P$ and such that $\mathcal B$ is at the left of the directed segment $\overrightarrow{AB}$. Let $F$ be the map that takes $A$ to $B$; clearly, $F$ is an orientation preserving homeomorphism of $\partial K.$

\begin{lemma}\label{Poncelet} The map $F$ of the curve $\partial K$ has an invariant measure and hence is conjugated to a circle rotation.
\end{lemma}

\emph{Proof.} Let $ds$ be the arc length element and let $\rho(A)=|AP|$ be the length of the tangent
segment from $A$ to the circle. Consider the measure
$d\mu=\frac{ds}{\rho}$, the measure of any Borel set $M\subset \partial K$ is given by $$\mu(M)=\int \limits _{M}\frac{ds}{\rho}.$$
The measure $\mu$ is an $F$-invariant measure, i.e., 
\begin{equation}\label{mu}
\mu(M)=\mu(F(M)).
\end{equation} 
Indeed, let $A'$ be a point infinitesimally close to $A$ and $B'=F(A').$
Then, by the condition that the product $|AP|\cdot |PB|$ is a constant
independent of the initial point $A$, the infinitesimal triangles
$\triangle AA'P$ and $\triangle BB'P$ are similar, and hence
$\frac{|AA'|}{|AP|}=\frac{|BB'|}{|B'P|}$, that is,
\begin{equation}\label{medida}
\frac{dA}{\rho(A)}=\frac{dB}{\rho(B)},
\end{equation} 
(we use the fact that the two
tangent segments to a circle from any point are equal). Integration of (\ref{medida}) gives (\ref{mu}). Now we apply Denjoy's theorem (see for instance Theorem 12.3 in \cite{Flatto}) and conclude that $F$ is conjugate to a circle rotation. \fin

\begin{figure}[H]
    \centering
    \includegraphics[width=.53\textwidth]{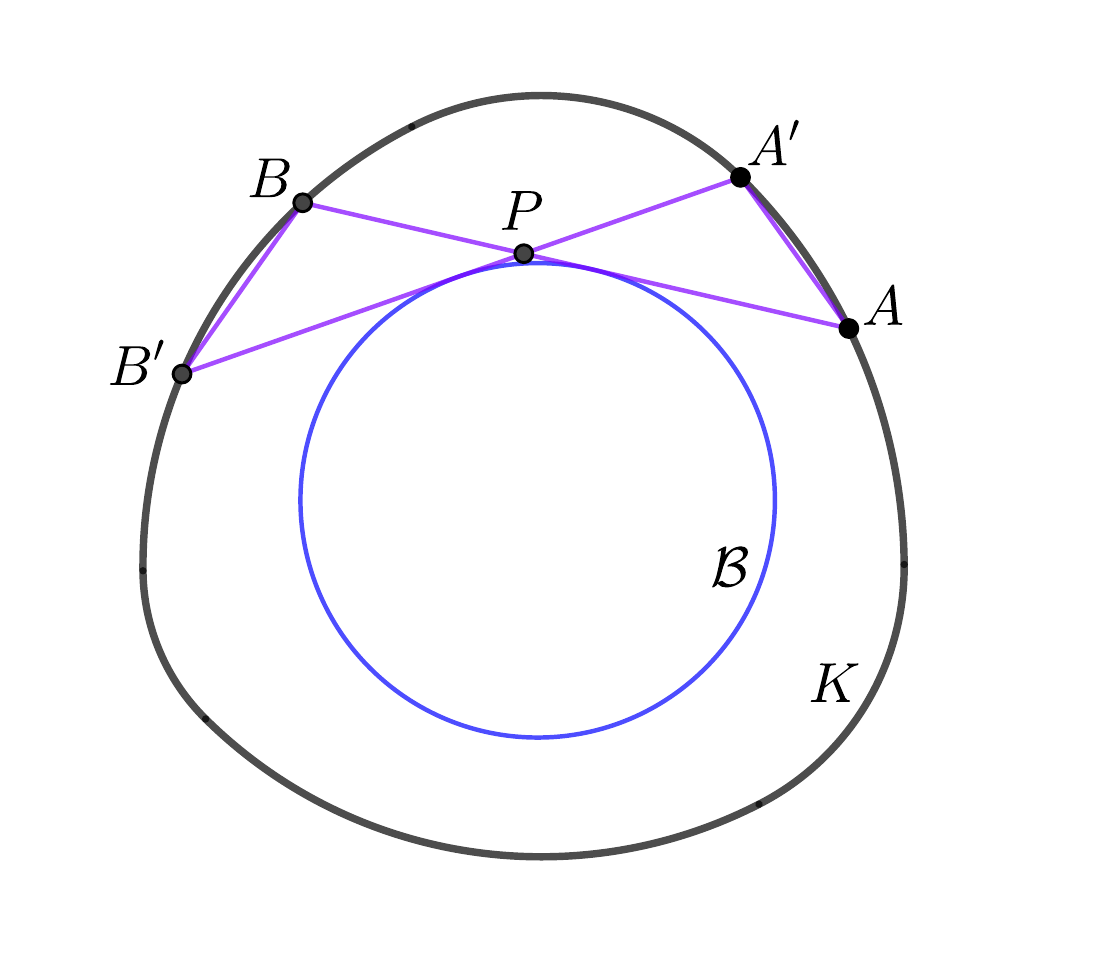}
    \caption{Triangles $\triangle AA'P$ and $\triangle BB'P$ are similar}
    \label{convexo}
\end{figure}

\begin{lemma}\label{rotacion}
Under the application of the map $F$ all
orbits are periodic with the same period or all are dense.
\end{lemma}

\emph{Proof.} This is a direct consequence of Lemma \ref{Poncelet}. \fin

We are ready now to prove the main result in this section.

\emph{Proof of Theorem \ref{principal}.} Let $A_0A_1$ be any chord of $K$ tangent to $\mathcal B$. Let $A_2$ be the point in $\partial K$ such that $A_1A_2$ is also tangent to $\mathcal B$, an so on. Denote by $\alpha=\measuredangle A_0OA_1$, if $\frac{\alpha}{\pi}$ is an irrational number, then the orbit of the point $A_0$ is dense in the boundary of $K$. Denote by $P_i$ the point of tangency between the segment $A_{i}A_{i+1}$ and $\mathcal B$, for every $i=0, 1, 2, \ldots$. Since $|P_0A_1|=|A_1P_1|$ and $|A_0P_0|\cdot |P_0A_1|=|A_1P_1|\cdot |P_1A_2|$, we have that $|A_0A_1|=|A_1A_2|$. Continuing in this way we see that all the chords $A_iA_{i+1}$ have the same length. By a theorem due to J. Barker and D. Larman given in \cite{Barker} we have that $K$ is a disc. Suppose now that $\frac{\alpha}{\pi}$ is a rational number, then the orbit of the point $A_0$ is periodic with period $n$, for some natural and fixed number $n\geq 3$, independently of the initial point $A_0$. Moreover, since the length of the chord $A_0A_1$ depends continuously on the point $A_0$, and every orbit has period $n$, we have that all the chords of $K$ and tangent to $\mathcal B$ are seen under the same angle $\alpha$. Now we will prove the following:

\textbf{Claim 1.} The length of all the chords of $K$ tangent to $\mathcal B$ is a constant.

\begin{figure}[H]
    \centering
    \includegraphics[width=.58\textwidth]{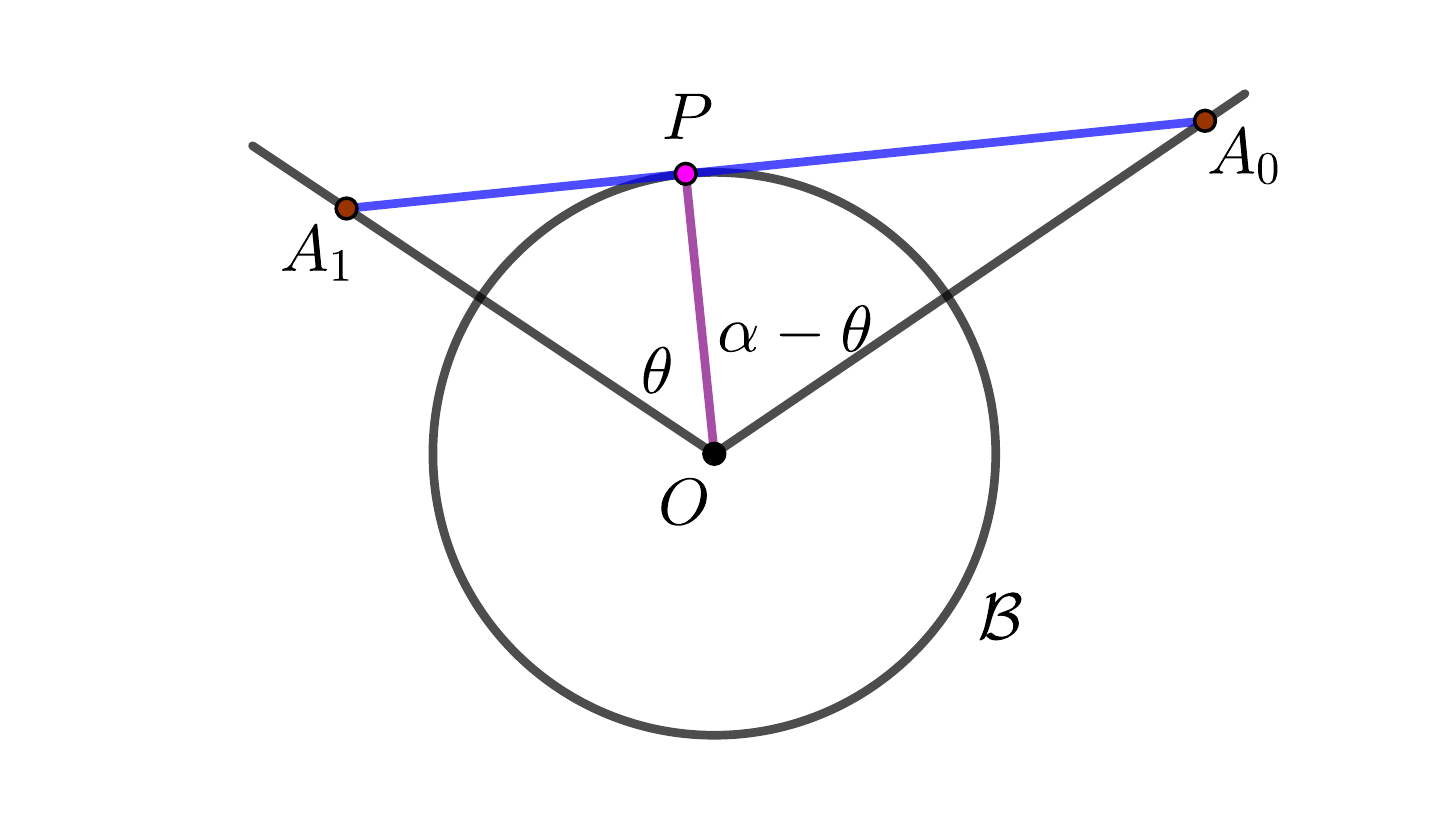}
    \caption{It holds that $|A_1P|\cdot|PA_0|=\tan \theta \cdot \tan (\alpha-\theta)$}
    \label{lema_cuerdas}
\end{figure}

\emph{Proof.} Suppose the point of contact between $A_0A_1$ and $\mathcal B$ is $P$ and the angle $\measuredangle A_1OP=\theta$. Without loss of generality we may suppose that the radius of $\mathcal B$ is 1 (see Fig. \ref{lema_cuerdas}). 

We have that $$|A_1P|\cdot|PA_0|=\tan \theta \cdot \tan (\alpha-\theta)=\lambda,$$ for some positive constant $\lambda$. We know that given the admissible value of $\lambda$, the following is true:

\textbf{Claim 2.} If $\alpha\neq \frac{\pi}{2}$ then there are only two values of $\theta$ that solves the equation $\tan \theta \cdot \tan (\alpha-\theta)=\lambda$, indeed, if one of the solutions is $\theta_0$ the other one is $\alpha-\theta_0$. 

\emph{Proof.} By the trigonometric identity for tangent of sum of angles we have that $$f(\theta):=\tan \theta \cdot \tan (\alpha-\theta)=1-\frac{\tan \theta + \tan(\alpha-\theta) }{\tan \alpha}.$$ Deriving the function $f$ we obtain $f'(\theta)=-\frac{1}{\tan\alpha}\left (\sec^2\theta -\sec^2(\alpha-\theta) \right ).$ If $\alpha=\pi/2$ then $\lambda$ must be equal to 1, $f'(\theta)=0$ and so every $\theta\in (0,\pi/4]$ gives an admissible solution. Now we consider the case $\alpha\neq \pi/2.$ Since the function $\sec^2(\theta)$ is an increasing function in the interval $[0,\pi/2)$, we have that $\sec^2\theta -\sec^2(\alpha-\theta)<0$, for every $\theta\in (0,\alpha/2],$ this implies that $f'(\theta)>0$, for every $\theta\in (0,\alpha/2].$ If follows that the function $f(\theta)$ is strictly increasing in the interval $\theta\in (0,\alpha/2],$ hence for every admissible value of $\lambda$ there are at most two values of $\theta\in(0,\alpha)$ which solves the equation. Indeed, there are exactly two values if $\alpha\neq \pi/2$, namely, $\theta$ and $\alpha-\theta.$ \fin

Claim 2 implies that there is only one value for the length of the chord $A_0A_1$. This ends the proof of Claim 1.\fin

Continuing with the proof of Theorem \ref{principal}, we have that in the case when $A_0$ has a periodic orbit,  all the chords of $K$ tangent to $\mathcal B$ have the same length. We use Barker-Larman's theorem and conclude that $K$ is a disc. 

Finally, if $\alpha=\frac{\pi}{2}$ then the orbit of $A_0$, for every $A_0\in \partial K$, has period equal to 4. This means that $O$ is the midpoint of $A_0A_2$, for every $A_0\in \partial K,$ i.e., $K$ has centre of symmetry at $O$. \fin

\begin{figure}[H]
    \centering
    \includegraphics[width=.55\textwidth]{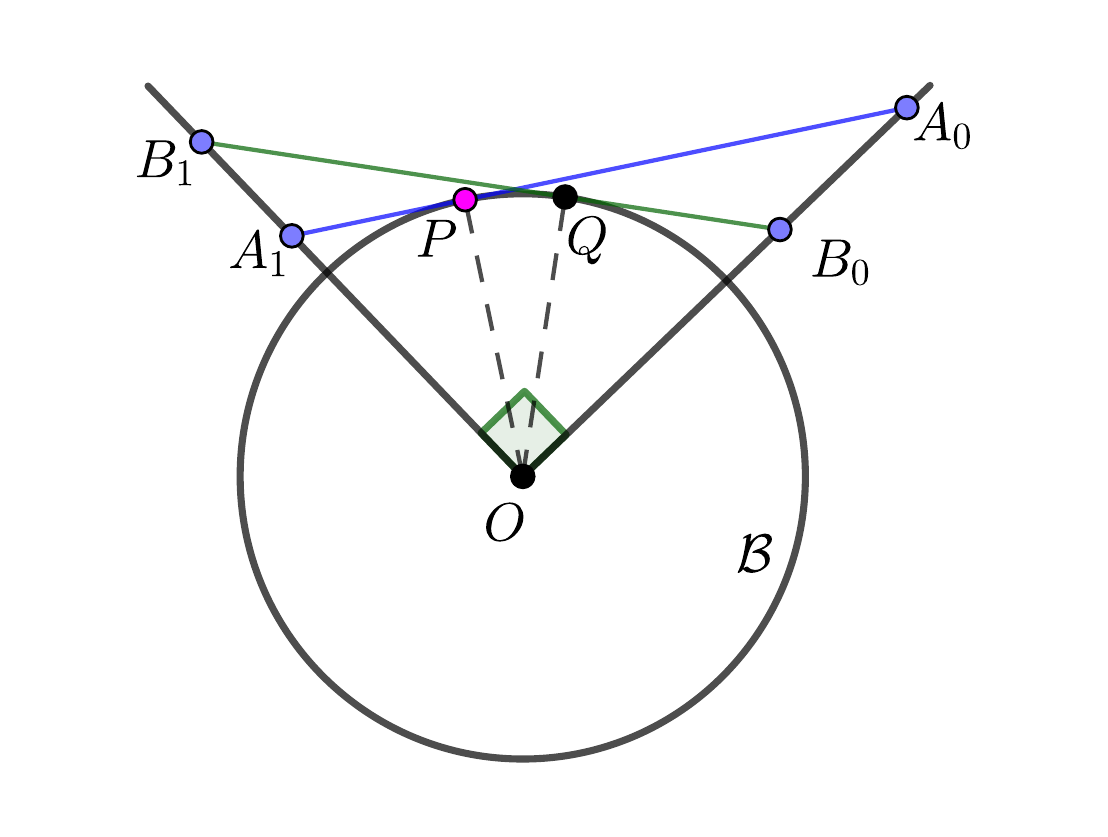}
    \caption{It holds that $|A_1P|\cdot|PA_0|=|B_1Q|\cdot|QB_0|$, however, $A_0A_1\neq B_0B_1$}
    \label{angulo_recto}
\end{figure}

\begin{remark}\label{ortoptica}
If $\alpha=\frac{\pi}{2}$ then $|A_1P|\cdot|PA_0|=|PO|^2$, but the same is true if $Q$ is any point in the arc of $\mathcal B$ contained in the angle $\measuredangle A_0OA_1$, i.e., we also have that $|B_1Q|\cdot|QB_0|=|QO|^2=|PO|^2$ (see Fig. \ref{angulo_recto}). This implies that the chords of $K$ tangent to $\mathcal B$ can have an infinite number of values for its length. Under this condition is not possible to say that $K$ is disc. 
\end{remark}

\begin{figure}[H]
    \centering
    \includegraphics[width=.58\textwidth]{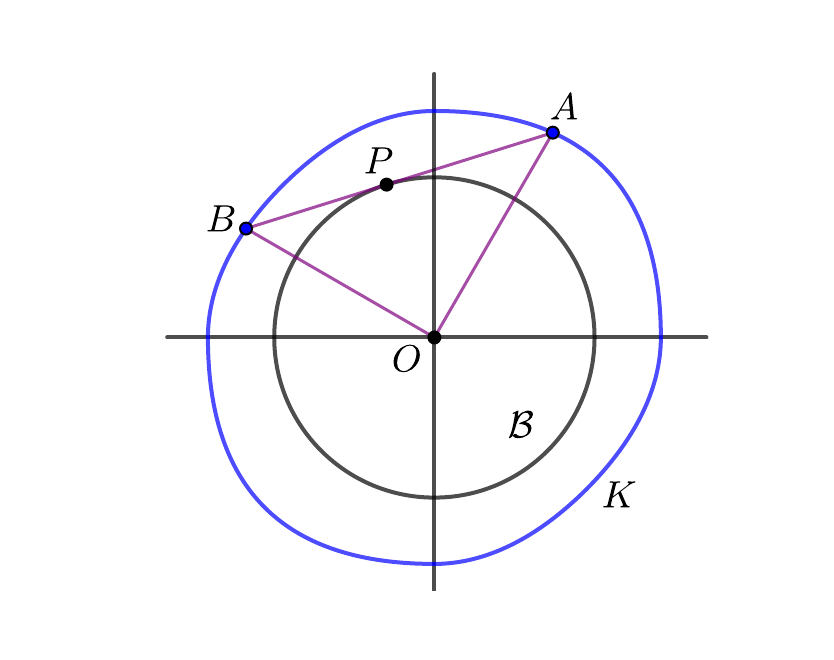}
    \caption{The body $K$ has an equipotential circle}
    \label{ejemplo}
\end{figure}

\begin{example}
Besides the example of the ellipse given in the introduction, we have the following: consider the arc of the parabola with focus in $(\frac{\sqrt 2}{2},\frac{\sqrt 2}{2})$ and directrix $y=-x+2\sqrt 2$, in the first quadrant of the Cartesian coordinate system. From every point $A$ in this arc, consider the tangent line to $\mathcal B$ and let $B$ be the point in this line such that $|BP|\cdot|PA|=1.$ The locus of the point $B$ is a convex curve in the second quadrant and shares the tangent line with the parabola at the point $(0,\sqrt 2)$. Now, consider the symmetric image of the union of these two arcs, with respect to $O$. In this way we obtain the boundary of a convex body $K$ which has $\mathcal B$ as its equipotential circle (see Fig. \ref{ejemplo}).
\end{example}

\begin{figure}[H]
    \centering
    \includegraphics[width=.7\textwidth]{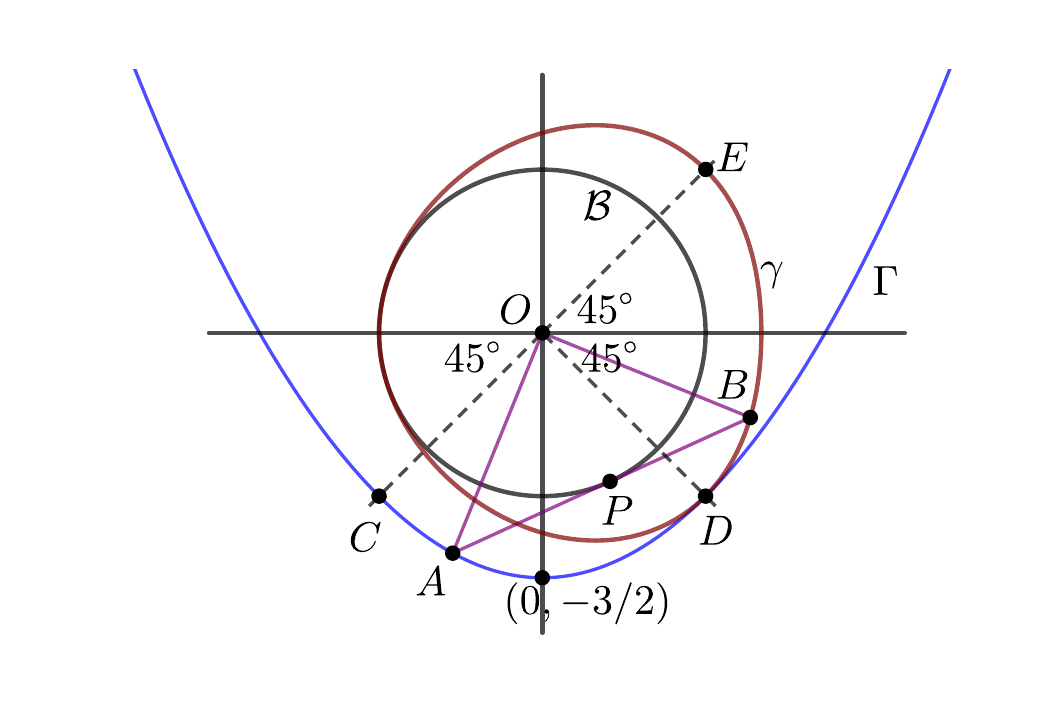}
    \caption{It holds that $|OB|=\frac{|OA|}{|AP|}$}
    \label{ejemplo2}
\end{figure}

We will explain why this body $K$ is convex. For the sake of clarity we rotate the parabola an angle $3\pi/4$ in the clockwise sense. Let $\Gamma$ be the parabola whose equation in Cartesian coordinates is $y=\frac{x^2}{2}-\frac{3}{2},$ and let $\gamma$ be the locus of the point $B$. The arc of the parabola in the first quadrant in Fig. \ref{ejemplo} corresponds to the  arc $\widehat{CD}$ in Fig. \ref{ejemplo2}. By the similarity of the triangles $\triangle OAP$ and $\triangle OBP$ we have that $$\frac{|OA|}{|OB|}=\frac{|AP|}{|OP|}=|AP|,$$ hence $|OB|=\frac{|OA|}{|AP|}.$ 

Let $\gamma'$ be the image of $\gamma$ under a rotation of an angle $\pi/2$ about $O$ in the clockwise sense. Let $B'$ be the image of the point $B$ under this rotation. We have that $|OB'|=\frac{|OA|}{|AP|}.$ The length of the segment $AP$ is $|AP|=\sqrt{-\frac{x^2}{2}+\frac{x^4}{4}+\frac{5}{4}}.$ If we parametrize $\Gamma$ and $\gamma'$ in polar coordinates by $(\theta, r(\theta))$ and $(\theta, r_1(\theta))$, respectively, we have that $$r_1(\theta)=\dfrac{r(\theta)}{\sqrt{-\dfrac{r^2(\theta)\cos^2\theta}{2}+\dfrac{r^4(\theta)\cos^4\theta}{4}+\dfrac{5}{4}}},$$ where $r(\theta)=\frac{\sin\theta+\sqrt{1+2\cos^2\theta}}{\cos^2\theta}.$ Since $r^2(\theta)\cos^2\theta -2r(\theta)\sin\theta -3=0,$ we obtain $$r_1(\theta)=\dfrac{r(\theta)}{\sqrt{r^2(\theta)\sin^2\theta +2r(\theta)\sin\theta +2}}.$$ We recall that the curvature in polar coordinates is calculated by $$\kappa(\theta)=\dfrac{\rho^2(\theta)-\rho(\theta)\rho''(\theta)+2(\rho'(\theta))^2}{((\rho'(\theta))^2 +\rho^2(\theta))^{\frac{3}{2}}},$$ so if we want to know that $\gamma'$ is a convex curve, we just need to see that\linebreak $r_1^2(\theta)-r_1(\theta)r_1''(\theta)+2(r'_1(\theta))^2>0$, for every $\theta\in[0,2\pi].$ This is indeed the case, we omit the simple but tedious computations, hence we have that $\gamma'$ is a convex curve, i.e., the boundary of a convex body. Moreover, $\gamma'$ pass through the points $C,$ and $D$ when $\theta=\frac{5\pi}{4},$ and $\theta=\frac{7\pi}{4},$ and share the tangent lines with $\Gamma$. This ensure us that the body $K$, constructed as we said before, is a convex body.

\begin{figure}[H]
    \centering
    \includegraphics[width=.84\textwidth]{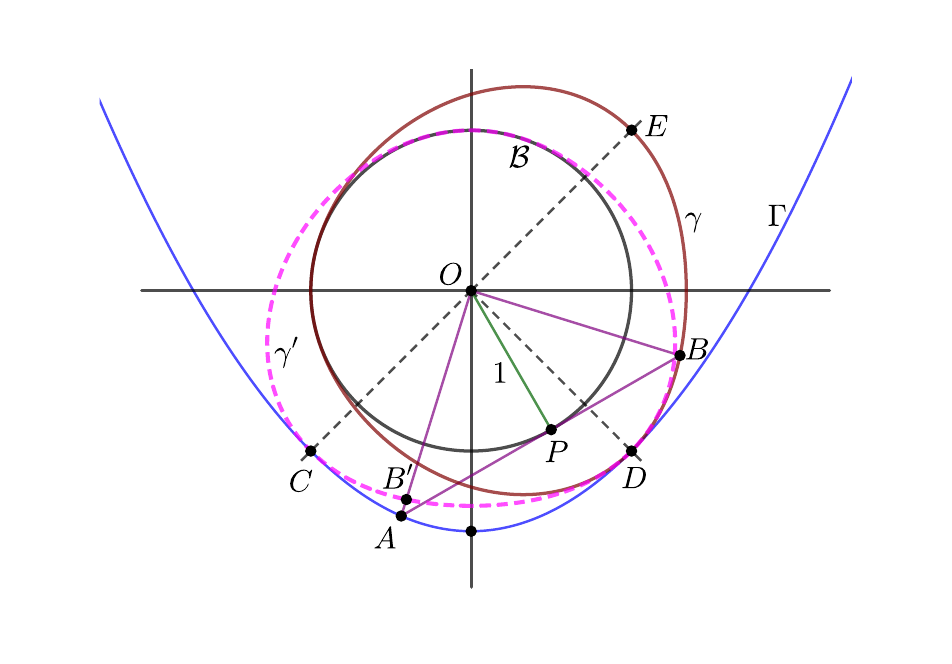}
    \caption{$\gamma$ is a convex curve}
    \label{ejemplo3}
\end{figure}

We prove that the only convex body with two equipotential circles is the disc.

\begin{theorem}\label{segundo}
Let $K\subset\mathbb R^2$ be a convex body which possesses two equipotential circles. Then $K$ is a disc.
\end{theorem}

\emph{Proof.} By Theorem \ref{principal} we have that $K$ is centrally symmetric and the two equipotential discs have centre at the centre of symmetry of $K$, to say $O$. By Lemma \ref{rotacion} and Remark \ref{ortoptica}, we have that for at least one of the equipotential discs it holds that the chords of $K$ and tangent to it, are seen from $O$ under a constant angle $\alpha\neq \frac{\pi}{2}$. By Theorem \ref{principal} we have that $K$ is a disc. \fin

\section{Some characterizations of the ellipsoid and the sphere}
In this section we derive from Theorem \ref{principal} some characterizations of the sphere and the ellipsoid.  

\begin{theorem}\label{falsopunto}
Let $K\subset\mathbb R^3$ be a convex body and let $P,$ $Q$ be two points in the interior of $K$. If every $2$-dimensional section of $K$ through $P$ or $Q$ possesses an equipotential circle, then $K$ is an ellipsoid.
\end{theorem} 
 
\emph{Proof.} By Theorem \ref{principal} we have that every section of $K$ through $P$ or $Q$ has a centre of symmetry. We recall the concept of a false centre: A convex body $L$ has a false centre $X$ if every section of $L$ through $X$ has centre of symmetry but $X$ is not centre of symmetry for $L$. Now we use the Theorem of the False centre of Aitchison, Petty, Rogers, Montejano, and Morales (see \cite{false_centre}, and \cite{Falso_MM}) which asserts that a convex body with a false centre must be an ellipsoid. Since either $P$ or $Q$ is not a centre of symmetry for $K$, we have that one of them is a false centre and then $K$ is an ellipsoid. \fin

The following result is a particular case of a theorem due to S. Olovjanishnikov \cite{Olovja}, however, we give our own proof for the case we need in some of the subsequent results. First we give some piece of notation: For $\upsilon\in \mathbb{S}^2$, we denote by $\upsilon^\perp$ the subspace  orthogonal to $\upsilon$ and by $E(\upsilon)$ the affine plane $\upsilon+\upsilon^\perp$. For $z\in \mathbb{R}^3 \setminus \mathbb{B}$, where $\mathbb B$ denote the unit ball of $\mathbb R^3$, we denote by $C_z$ the union of all the supporting lines of $\mathbb{S}^2$ passing through $z$.

\begin{theorem}\label{otravezpami}
Let $K \subset \mathbb{R}^3$ be a convex body. Suppose that 
$ \mathbb{S}^2 \subset \emph{int}\, K$ and, for every $\upsilon\in \mathbb{S}^2$, the section $E(\upsilon) \cap K$ is centrally symmetric with centre at $\upsilon$. Then $K$ is a ball concentric with $\mathbb{S}^2$.
\end{theorem}

\begin{figure}[H]
    \centering
    \includegraphics[width=.78\textwidth]{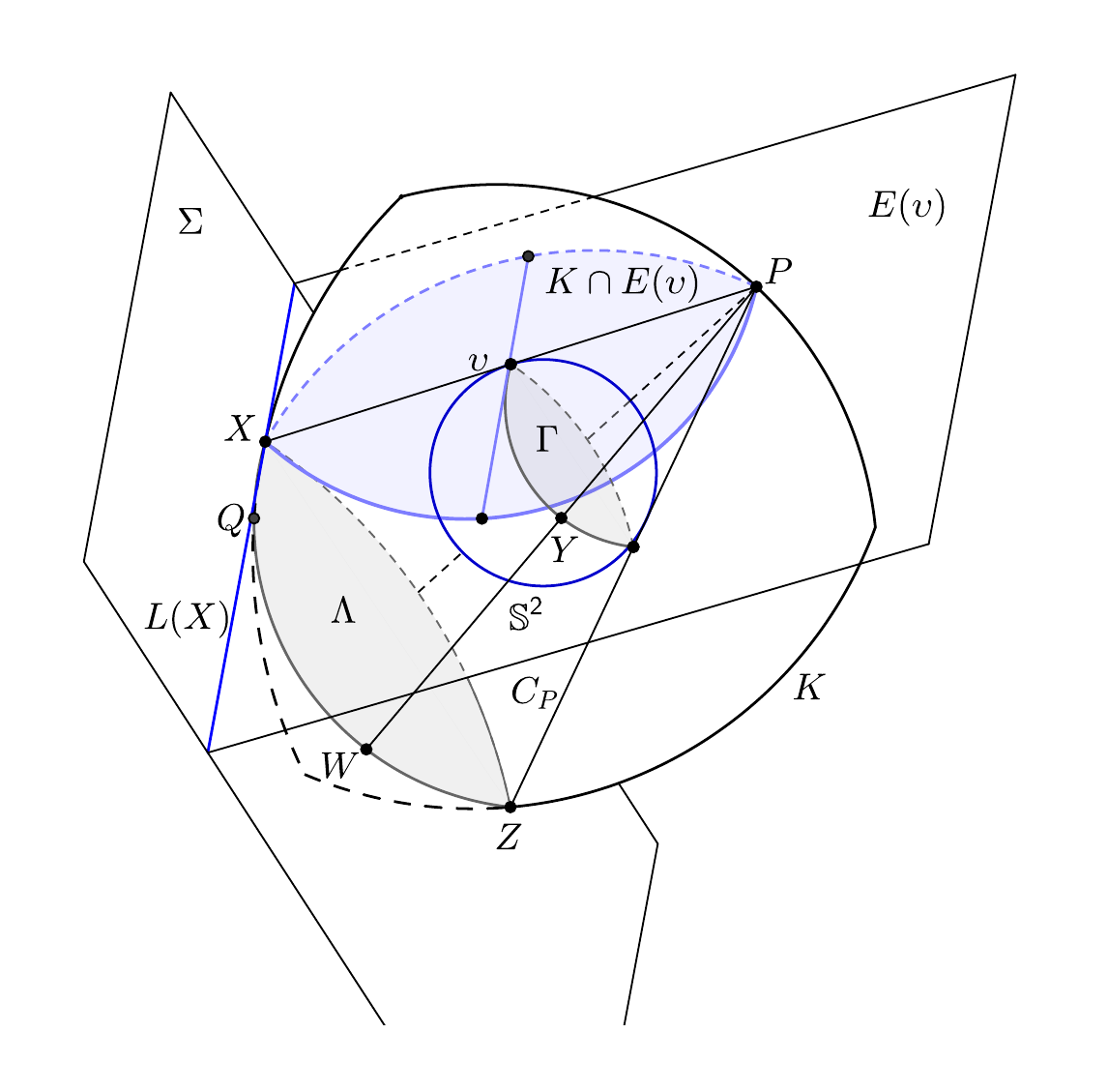}
    \caption{The section $E(\upsilon) \cap K$ is a circle}
    \label{lema_Olovja}
\end{figure}

\emph{Proof.} First of all, we are going to prove that, for all $P \in \partial K$, the projection of $\mathbb{S}^2$ from $P$ onto $\partial K$ is a circle. Let $P \in \partial K$ and let $\Gamma$ be the circle $C_{P} \cap \mathbb{S}^2$ (see Fig. \ref{lema_Olovja}). We will show that the set $\partial K \cap [C_{P}  \setminus \{ P \}]$, which is denoted by 
$\Lambda$, is a circle. Let $W\in \Lambda$. Since the segment $WP \subset C_{P}$, the line $WP$ makes contact with $\mathbb{S}^2$ in a point, say $Y$. It is clear that $Y\in \Gamma$. In virtue that $E(Y) \cap K$ has centre at $Y$, $Y$ must be the midpoint of the chord $PW$. It follows that 
$\Lambda =\phi (\Gamma)$, where $\phi: \mathbb{R}^3 \rightarrow \mathbb{R}^3$ is an homothety with centre $P$ and ratio of homothety equal to 2. Therefore, $\Lambda$ is a circle and is contained in a plane parallel to the plane of  $\Gamma$.

Now we are going to prove that, for every $\upsilon\in \mathbb{S}^2$, the 
section $E(\upsilon) \cap K$ has the following property: \textit{for every point $X \in \partial [E(\upsilon) \cap K]$ there exist a supporting line $L(X)$ of $E(\upsilon) \cap K$ such that $L(X)$ and the line $X\upsilon$ are orthogonal}. It is well known that such property characterizes the disc \cite{Toponogov}.

Consider an arbitrary point $\upsilon\in \mathbb{S}^2$ and a point $X \in \partial [E(\upsilon) \cap K]$. We denote by $P$ the intersection of ray $\overrightarrow{X\upsilon}$ with $\partial K$. By the first part of the proof, there exist a plane $\Sigma$ such that $\partial K \cap [C_{P}  \setminus \{ P \}]= \Sigma \cap \partial K$ is a circle, say again $\Lambda$, and $X\in \Sigma \cap \partial K$. Since $E(\upsilon)$ is a supporting plane of $\mathbb {S}^2$ (and of $C_{P}$), the line $E(\upsilon) \cap \Sigma$ is a supporting line of $\Lambda$. The line $E(\upsilon) \cap \Sigma$ must be also a supporting line of $E(\upsilon) \cap K$, otherwise, there is a point $Q\in  E(\upsilon) \cap K$ such that $Q\in E(\upsilon) \cap \Sigma$ and $X\neq Q$. Since $Q\in  E(\upsilon) \cap \Sigma$ it follows that $Q\in \Sigma$. We have that the points of $K$ which are in $\Sigma$ are those points in $\Lambda$, we conclude that $Q\in \Lambda$, i.e., there is a point $Q$ in the supporting line $E(\upsilon) \cap \Sigma$ of the circle $\Lambda$ different than $X$, which is absurd. Hence the line $E(\upsilon) \cap \Sigma$ is a supporting line of $E(\upsilon) \cap K$. Finally, we observe that, since $C_{P}$ is a right circular cone, then $E(\upsilon) \cap \Sigma$ is orthogonal to the line $X\upsilon$. 

The proof of the theorem will be complete as soon as we prove the following.

\textbf{Claim 3.} If all the sections of a convex body $K\subset \mathbb R^3$, given by supporting planes of a sphere $S\subset \text{int} K$, are circles and the centre of each section is the contact point of the supporting plane with the sphere, then $K$ is a ball concentric with $S$. 

\emph{Proof.} Suppose the centre and radius of $S$ are $O$ and $r$, respectively. Let $H$ be an arbitrary plane tangent to $S$ at the point $A$ and denote by $S_A$ the circle $H\cap K$. Let $K_A$ be the open region of $\partial K$ which is separated of $S$ by mean of the plane $H$. Consider any point $X\in K_A$ and let $\Gamma_X$ be the any section of $K$ tangent to $S$ and passing through $X$. Clearly, $\Gamma_X$ intersects the circle $S_A$ in two points, to say $B,$ and $C$ and is tangent to $S$ in a point $D$. The radius of $S_A$ and $S_D$ is the same and is equal to $|OB|^2-r^2$, hence, $\overline{K_A}$ is a closed cap of a sphere with centre $O$ and radius $|OB|$. Let $Z$ be any point in $S$ such that $K_Z\cap K_A\neq \emptyset$. Again, $\overline{K_Z}$ is a closed cap of the sphere with centre $O$ and radius $|OB|$. Continuing in this way, we obtain a closed cover of $\partial K$, since $\partial K$ is a compact set, we have that there exists a finite subcover of $\partial K$. Every region of this subcover is a cap of the sphere with centre $O$ and radius $|OB|$, therefore, $K$ is a ball concentric with $S$. \fin

The proof of Theorem \ref{otravezpami} is now complete. \fin

As a corollary of Theorems \ref{principal} and \ref{otravezpami} we have the following.

\begin{theorem}\label{esfera}
Let $K\subset\mathbb R^3$ be a convex body which contains a ball $\mathcal B$ in its interior. Suppose every section of $K$, and tangent to $\mathcal B$, has the contact point as the centre of an equipotential circle. Then $K$ is a ball.
\end{theorem}

\emph{Proof.} By Theorem \ref{principal} we have that every section of $K$ tangent to $\mathcal B$ has the point of contact as centre of symmetry. By Theorem \ref{otravezpami}, we have that $K$ and $\mathcal B$ are concentric balls, therefore, $K$ is a ball. \fin

\begin{remark}
If we consider any convex body $L$ in Theorem \ref{esfera} instead of the ball $\mathcal B$, the conclusion is that $K$ and $L$ are homothetic ellipsoids. This conclusion follows if we apply Olovjanishnikov's theorem (see \cite{Olovja}) which asserts: \emph{Let $K$ and $L$ be convex bodies in $\mathbb R^n$, with $L\subset\text{int}\, K$, such that every hyper-section of $K$ tangent to $L$ is centrally symmetric and its centre belongs to $\partial L$, then $K$ and $L$ are homothetic ellipsoids. } \end{remark}

A mixture between Theorems \ref{falsopunto} and \ref{esfera} is the following.

\begin{theorem}\label{esfera2}
Let $K\subset\mathbb R^3$ be a convex body which contains a ball $\mathcal B$ in its interior and let $P$ be a point in the interior of $\mathcal B$, different from the centre of $\mathcal B$. Suppose every section of $K$ through $P$ has as an equipotential circle, the intersection of the section with $\partial \mathcal B$. Then $K$ is a ball.
\end{theorem}

\emph{Proof.} First of all, we are going to prove that $P$ is a false centre of $K$. By Theorem \ref{principal} we have that all sections of $K$ passing through $P$ are centrally symmetric. We denote by $O$  the centre of $\mathcal B$. Since $P\neq O$, there exist a line $L$ through $O$ such that $P$ does not belong  to $L$. Let $\Pi$ be the plane orthogonal to $L$ through $P$ and denote by $T$ the intersection $\Pi \cap L$. For one hand, since $P$ does not belong to $L$, $P\neq T$, on the other hand, $T$ is the centre of the section $\Pi \cap \mathcal B$ and, consequently, is the centre of $\Pi \cap K$. Hence $P$ is a false centre of $K$. By the the False centre Theorem (see \cite{false_centre}, and \cite{Falso_MM}) $K$ is an ellipsoid. We denote the centre of $K$ by $R$. If $O\neq R$, then the sections $\text{aff}\{P,O,R\}\cap \mathcal B$ and $\text{aff}\{P,O,R\}\cap K$ (where $\text{aff}\{P,O,R\}$ denotes the plane through $P,$ $O,$ and $R$), would have centres at $O$ and at $R$, respectively. This would contradict Theorem \ref{principal} (if $P,O, R$ were collinear we would consider any plane containing the line determined by this points).
 
On the other hand, is well known that the locus $\Omega_E(P)$ of the mid points of the chords of the ellipsoid $E$, with centre at $O$, passing through the point $P\in \text{int} E$ is an ellipsoid homothetic to $E$ and $P,O\in \Omega_E(P)$. We denote by  $\Delta_E(P)$ the locus of the centres of the sections of $E$ passing through $P$. It is easy to see that 
$\Omega_E(P)=\Delta_E(P)$. Thus $\Delta_E(P)$ is an ellipsoid homothetic to $E$ and $P,O\in \Delta_E(P)$. By Theorem \ref{principal} it follows that $\Delta_K(P)=\Delta_S(P)$, that is, $\Delta_K(P)$ is a sphere. Since $\Delta_K(P)$ and $K$ are homothetic, we conclude that $K$ is a ball. \fin

For the next result we need the following definition: We say that a sphere $\mathcal S$ in the interior of a convex body $K\subset\mathbb R^3$, is an \emph{equipotential sphere} if for every point $P\in\mathcal S$, the tangent plane at $P$ cuts a section $K(P)$ of $K$, such that $P$ is an equipotential point of $K(P)$ with an equipotential constant $\lambda (P)$ (depending on the point $P$).

\begin{theorem}\label{esfera_equipotencial}
Let $K\subset\mathbb R^3$ be a convex body which possesses and equipotential sphere $\mathcal S$ with $\lambda (P)=\lambda_0$, for every $P\in \mathcal S$. Then $K$ is a ball.
\end{theorem}

\emph{Proof.} Let $P$ be any point in the interior of $\mathcal S$, different from the centre of $S$. Consider any section $\mathcal H$ of $K$ through $P$, since $\mathcal S$ is an equipotential sphere with an equipotential constant independent of the point in $\mathcal S$, we have that $\mathcal H$ has the circle $\mathcal H\cap \mathcal S$ as an equipotential circle. We have the conditions of Theorem \ref{esfera2}, therefore, $K$ must be a ball. \fin
 
\section{Convex bodies with equireciprocal circles}
Let $\mathcal E$ be an ellipse with foci $P$ and $Q$. Consider any chord $AB$ through $P$ (or $Q$), then the value of $\frac{1}{|AP|}+\frac{1}{|PB|}$ is equal to a constant $\lambda$. If a point $P$ for a convex figure $K$ has the property that for any chord $AB$ through $P$ the value of $\frac{1}{|AP|}+\frac{1}{|PB|}$ is equal to a constant, we then say that $P$ is an \emph{equireciprocal point}. In \cite{Falconer}, K. Falconer proved that a $C^2$ convex curve with two equireciprocal points must be an ellipse. However, if the differentiability hypothesis is omitted, then there exist some convex curves, different from ellipses, which have two equireciprocal points. In this section we extent the notion of equireciprocal point: we say that a circle $\mathcal B$ in the interior of a convex body $K\subset \mathbb R^2$ is an \emph{equireciprocal circle} if for every chord $AB$ of $K$, tangent to $\mathcal B$ at a point $P$ (depending on the chord $AB$), it holds that $\frac{1}{|AP|}+\frac{1}{|PB|}=\lambda$, for some constant $\lambda$. It is natural to ask about the existence of a convex body $K$, different from the disc, which has an equireciprocal disc. The answer is given in the following.

\begin{theorem}\label{equireciproco}
Let $K\subset\mathbb R^2$ be a convex body containing an equireciprocal disc $\mathcal B$ with centre $O$ in its interior. Then $K$ is a disc centred at $O$. 
\end{theorem}
  
\emph{Proof.} The proof is completely analogous to the proof of Theorem \ref{principal} with only one change: in Claim 2 we use the equation $\cot \theta +\cot(\alpha-\theta)=\lambda$, instead of the equation $\tan\theta\cdot \tan(\alpha-\theta)=\lambda.$ In this case, the derivative of $f$ is $$f'(\theta)=-\csc^2\theta +\csc^2(\alpha-\theta),$$ which is negative for $\theta$ in the interval $(0,\alpha/2)$, independently of the value of $\alpha\in(0,\pi)$. Hence, for the admissible value $\lambda$ there is only one solution of the mentioned equation. In other words, for every value of $\alpha$, there is only one possible length for the chord $AB$. Again, by the theorem of Barker and Larman \cite{Barker}, we conclude that $K$ is a disc centred at $O$. \fin  

We say that a sphere $\mathcal S$ in the interior of a convex body $K\subset\mathbb R^3$, is an \emph{equireciprocal sphere} if for every point $P\in\mathcal S$, the tangent plane at $P$ cuts a section $K(P)$ of $K$, such that $P$ is an equireciprocal point of $K(P)$ with an equireciprocal constant $\lambda (P)$ (depending on the point $P$). As a first corollary of Theorem \ref{equireciproco} we have the following.

\begin{theorem}\label{esfera3}
Let $K\subset\mathbb R^3$ be a convex body which contains an equireciprocal sphere $\mathcal S$ in its interior, with $\lambda(P)=\lambda_0$ for every $P\in\mathcal S$. Then $K$ is a ball.
\end{theorem}

\emph{Proof.} Every section $\mathcal H$ of $K$ which intersects the interior of $\mathcal S$, has $\mathcal H\cap \mathcal S$ as an equireciprocal circle, hence $\mathcal H$ is a disc with centre at the centre of $\mathcal H\cap \mathcal S$. Now, let $P$ be any point in the interior of $\mathcal S$. We know that there exists an affine diameter of $K$ through $P$ (see for instance \cite{Soltan}), to say, $AB$. Every section $\mathcal H$ of $K$ which contains to $AB$, is a disc with centre at the midpoint of $AB$ and radius equal to $|AB|/2$. It follows that $K$ is a ball with centre at the midpoint of $AB$ and radius $|AB|/2$. \fin  
  
For the following results we need a new definition (see \cite{JMM}): Let $\delta :\mathbb{S}^2\longrightarrow\mathbb{R}$ be a continuous function such that $\delta (-x)=-\delta (x)$. We denote by $\Sigma$ the following set of planes in $\mathbb{R}^3$:
$$\Sigma =\{ H_y = \{x\in\mathbb{R}^3|\langle x,y\rangle = \delta (y)\}\}_{y\in\mathbb{S}^2},$$ and we say that $\Sigma$ a $2$-cycle of planes. If for a given convex body $K$, it holds that $H_p\cap H_q\cap \text{int}\, K\neq\emptyset$ for every two planes $H_p,H_q\in\Sigma$, we say that $\Sigma$ is a $2$-cycle of planes for $K$. In \cite{Larman}, D. Larman, L. Montejano, and E. Morales proved some results which characterizes the ellipsoid in terms of the properties of the sections, cut by the members of a $2$-cycle of planes, in a given convex body. 

We have the following result which is interesting by itself.

\begin{theorem}\label{2ciclo}
Let $K\subset\mathbb{R}^{3}$ a convex body and let $\Sigma$ a $2$-cycle of planes for $K$. Suppose that for each $\Gamma\in\Sigma$ we have that $\Gamma\cap K$ is a disc. Then $K$ is a ball.
\end{theorem}

\emph{Proof.} Consider a unitary vector $u\in\mathbb{S}^2$, we may consider that $\delta(u)$ is the distance with sign from $H_u$ to the origin $O$. Recall that $H_u$ is the plane in $\Sigma$ that is orthogonal to $u$. Let $L$ be an affine diameter of $K$ parallel to $u$. 

\textbf{Claim 4.} There exists $\Pi\in\Sigma$ such that $L\subset\Pi$.

\emph{Proof.} Let $\Sigma_u$ be the subset of $\Sigma$ such that for every $H_y\in\Sigma_u$ it holds that $y$ is orthogonal to $u$. Denote by $S(u)=\mathbb S^2\cap u^{\bot}.$ Let $\delta_u$ be the restriction of $\delta$ to the set $S(u)$. Let $O'$ be the midpoint of the affine diameter $L$, and let $\sigma: S(u)\longrightarrow\mathbb R$ be the function such that $\sigma (x)$ is the distance with sign from $O'$ to $H_x$. Since $\delta_u$ is a continuous function we have that $\sigma$ is also a continuous function. Moreover, $\sigma(-x)=-\sigma(x)$, i.e., $\sigma$ is an odd function. By the Borsuk-Ulam theorem (see for instance \cite{Matousek}), there exists $x_0\in S(u)$ such that $\sigma(x_0)=0$, this means that the plane $H_{x_0}$ contains to $L$. \fin

We will show that there exists a plane of symmetry of $K$, orthogonal to $u$. 
Let $\Delta$ be the plane orthogonal to $u$ that passes through the midpoint of $L$. We assert that $\Delta$ is the plane of symmetry we are looking for. For $X\in bd(K-\Pi)$, denote as $L_X$ the chord of $K$ through $X$ parallel to $u$. Again, since $\Sigma$ is a $2$-cycle, there exists $\Pi_X\in\Sigma$ such that $L_X\subset\Pi_X$. 

Since the line $\Pi\cap\Pi_X$ has inner points of $K$, then the discs $\Pi\cap K$ and $\Pi_X\cap K$ have the common chord $(\Pi\cap\Pi_X)\cap K$. The plane $\Delta$ pass through the midpoint of the chord $L_X\subset\Pi_X\cap K$, also, the centre of $\Pi_X\cap K$ is contained in $\Delta$. Hence, the chord $\Delta\cap\Pi_X$ is a line of symmetry of the disc $\Pi_X\cap K$, and the second point of intersection of $L_X$ with $\Pi_X\cap K$, say $Y$, is at the same distance to the line $\Delta\cap\Pi_X$ as $X$. It follows that $\Delta$ is a plane of symmetry of $K$. We have proved that for every affine diameter $L$ of $K$, the plane orthogonal to $L$ through its midpoint is a plane of orthogonal symmetry of $K$, therefore,  $K$ must be a ball. \fin
  
As a corollary of Theorems \ref{equireciproco} and \ref{2ciclo} we have our last result. 
  
\begin{corollary}
Let $K\subset\mathbb{R}^{3}$ a convex body and let $\Sigma$ a $2$-cycle of planes. Suppose that for each $\Gamma\in\Sigma$ we have that $\Gamma\cap K$ possesses an equireciprocal circle. Then $K$ is a ball.
\end{corollary}
  
\emph{Proof.} The proof is just an application of Theorems \ref{equireciproco} and \ref{2ciclo}. \fin
  
\section{Further comments and problems}
In this section we give some problems and final comments.

\textbf{Problem 1.} Let $K\subset\mathbb R^2$ be a convex body which has a disc $\mathcal B$ in its interior. Suppose every chord of $K$ which is tangent to $\mathcal B$ subtends a constant angle $\alpha\neq \pi/2$ from the centre $O$ of $\mathcal B$. Is it $K$ a disc or an ellipse?

It is known that if $\alpha$ is a fixed angle, $\Omega$ an ellipse with foci $A$ and $B$, then the envelope of the chords of $\Omega$, which are seen under a constant angle $\alpha$ from a focus, to say $A$, is another ellipse with one focus in $A$ (see \cite{Smith}).

\begin{figure}[H]
    \centering
    \includegraphics[width=.68\textwidth]{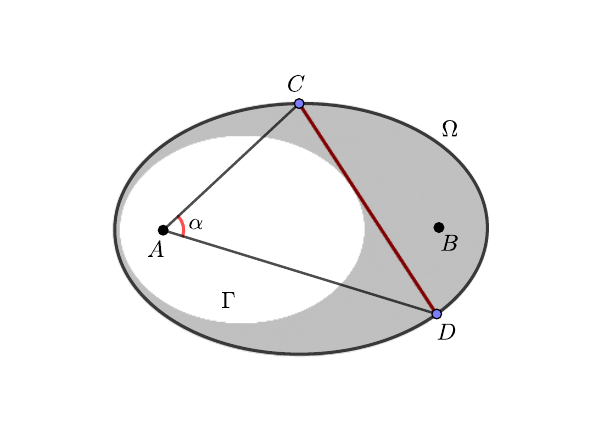}
    \caption{The envelope of the chords $CD$ is an ellipse}
    \label{elipse_cuerda}
\end{figure}

\textbf{Problem 2.} Let $K\subset\mathbb R^3$ be a convex body which has a ball $\mathcal B$ in its interior. Suppose every section of $K$ which is tangent to $\mathcal B$ subtends a solid constant angle $\alpha$ from the centre $O$ of $\mathcal B$. Is it $K$ a ball?

\textbf{Problem 3.} Is there an $O$-symmetric convex body $K\subset\mathbb R^3$, different from a ball,  which contains a ball $\mathcal B$ centred at $O$, such that every section of $K$ through $O$ has as an equipotential circle, the intersection of the section with $\partial \mathcal B$?

\textbf{Problem 4.} Let $K\subset\mathbb R^3$ be a convex body which possesses and equipotential sphere $\mathcal S$, with the equipotential constant $\lambda(P)$ not necessarily equal to a fixed number $\lambda_0$. Is it $K$ a ball?

\textbf{Problem 5.} Let $K\subset\mathbb R^3$ be a convex body which possesses and equireciprocal sphere $\mathcal S$, with the equireciprocal constant $\lambda(P)$ not necessarily equal to a fixed number $\lambda_0$. Is it $K$ a ball?

In a second paper on this topic, we will give more characterizations of the sphere and ellipsoid in spaces of higher dimension. In these characterizations we will use appropriate definitions for equipotential and equireciprocal spheres. 

\section{Acknowledgment}
We thank Diana Verdusco, Alejandro Estrada, Ulises Velasco, Carlos Yee, Luguis de los Santos,  and Sergei Tabachnikov, for many helpful discussions about the topic of this work. 

\end{document}